\documentclass[12pt]{article}

\usepackage[dvips]{graphicx}
\usepackage{float,latexsym,amsmath,amsthm,amsfonts,natbib}
\usepackage{amsthm}
\usepackage{color}

\newtheorem{theo}{Theorem}[section]
\newtheorem{lemma}{Lemma}[section]
\newtheorem{prop}{Proposition}[section]

{
\theoremstyle{definition}
\newtheorem{remark}{Remark}[section]
\newtheorem{defi}{Definition}[section]
}

\newcommand{\E}{\mathbb{E}}
\newcommand{\V}{\mathbb{V}}

\newcommand{\Es}{\mathbb{E}_{\star}}
\newcommand{\Vs}{\mathbb{V}_{\star}}

\newcommand{\PP}{\mathbb{P}}

\newif\ifnotes
\notestrue
\newcounter{mnotecount}[section]

\def\H{\mbox{$\mathcal H$}}
\DeclareMathOperator*{\argmin}{argmin}
\DeclareMathOperator*{\C}{Cov}

\date{}   

\begin{document}

\begin{center}

{\sc \Large Multiplicative Bias Corrected Nonparametric Smoothers}
%with Application to Nuclear Energy Spectrum Estimation}

\vspace{0.7cm}

%N. Hengartner$^{\mbox{\footnotesize a}}$, E. Matzner-L\o ber$^{\mbox{\footnotesize b}}$, L. Rouvi\`ere$^{\mbox{\footnotesize c,}}\footnote{Corresponding author: laurent.rouviere@ensai.fr}$ \footnote{Research partially supported by the French "Agence Nationale pour la Recherche" under grant ANR-09-BLAN-0051-02 "CLARA"} and T. Burr$^{\mbox{\footnotesize a}}$

\vspace{0.5cm}

%$^{\mbox{\footnotesize a}}$ Los Alamos National Laboratory\\ NM 87545, USA
\bigskip

%$^{\mbox{\footnotesize b}}$ Universit\'e de Rennes \\
% 35043 Rennes Cedex, France\\
\bigskip

%$^{\mbox{\footnotesize c}}$ CREST (Ensai), IRMAR, UEB 
%\\
%Campus de Ker Lann\\
%Rue Blaise Pascal - BP 37203\\
%35172 Bruz Cedex, France\\
\smallskip

\end{center}
\begin{abstract}
The paper presents a multiplicative bias reduction estimator for nonparametric regression. The approach consists to apply a multiplicative bias correction to an oversmooth pilot estimator. In \cite{burr2010}, this method has been tested to estimate energy spectra. For such data set, it was observed that the method allows to decrease bias with negligible increase in variance. In this paper, we study the asymptotic properties of the resulting estimate and prove that this estimate has zero asymptotic bias and the same asymptotic variance as the local linear estimate. Simulations show that our asymptotic results are available for modest sample sizes. 
%We also illustrate the benefit of this new method on nuclear energy spectrum estimation.
\medskip

\noindent{\bf Index terms: } Nonparametric regression, bias reduction, local linear estimate.
\end{abstract}

%\maketitle

\section{Introduction}
In nonparametric regression, the bias-variance tradeoff of linear smoothers such as kernel-based regression smoothers, wavelet based smoother or spline smoothers, is generally governed by a user-supplied parameter. This parameter is often called the bandwidth. As an example, assuming that the regression function $m$ is twice continuously differentiable, the local linear smoother with bandwidth $h$ and kernel $K$ has conditional bias
$$\frac{h^2}{2}m''(x)\int u^2K(u)\,\mathrm{d}u+\mathrm{o}_p(h^2)$$
and conditional variance 
$$\frac{1}{nh}\frac{\sigma^2(x)}{f(x)}\int K^2(u)\,\mathrm{d}u+\mathrm{o}_p\left(\frac{1}{nh}\right)$$
where $f$ stands for the density of the (one-dimensional) explanatory variable $X$ and $\sigma^2(x)$ is the conditional variance of the response variable given $X=x$. See for example the book of \cite{fan+1996}. Since the bias increases with the second order derivative of the regression function, the local linear smoother tends to under-estimate in the peaks and over-estimate in the valleys of the regression function. See for example \cite{simonoff1996,wand+1995,scott1992}. This behavior in peaks and valleys causes some trouble in some practical applications, such as for the estimation of energy spectrum.

\bigskip
The decay of radioactive isotopes often generates gamma particles whose energy can be measured using specialized detectors. Typically, these detectors count the number of particles in various energy bins over short time intervals such as one to ten minutes. This enables estimation of the energy distribution of the emitted particles, which is called the energy spectrum.  For low or medium resolution detectors, the spectrum is typically composed of multiple broad peaks whose location and area characterize the radio-isotope. That is we not only want the locations of the peaks, but also the shape and amplitudes in peak regions.
\medskip

Because the actual bin counts are noisy, and the energy spectrum is fairly smooth, it has been proposed to estimate the energy spectrum using nonparametric smoothing techniques (\cite{sullivan++2006,gang+++2004}). However, 
%it is known that many  classical smoothers, such as kernel-based regression smoothers, $k$-nearest neighbors, and smoothing splines, typically under-estimate in the peaks and over-estimate in the valleys of the regression function.  See for example \cite{simonoff1996,fan+1996,wand+1995,scott1992}.
the bias of the nonparametric smoothers degrades isotope identification performance for any algorithm that includes peak area or ratios of areas (\cite{casson+2006}) and motivates studying methods to reduce bias.
\medskip

All nonparametric smoothing methods are generally biased. There are many approaches to reducing the bias, but most of them do so at the cost of an increase in the variance of the estimator.  For example, one may chose to undersmooth the energy spectrum. Undersmoothing will reduce the bias but will have a tendency of generating spurious peaks.  One can also use higher order smoothers, such as local polynomial smoother with a polynomial of order larger than one.  While again this will lead to a smaller bias, the smoother will have a larger variance. Another approach is to start with a pilot smoother and to estimate its bias by smoothing the residuals (\cite{marzio+2007}). Subtracting the estimated bias from the smoother produces a regression smoother with smaller bias and larger variance. For the estimation of an energy spectrum, the additive bias correction and the higher order smoothers have the unfortunate side effect of possibly generating a non-positive estimate.
\medskip

An attractive alternative to the linear bias correction is the multiplicative bias correction pioneered by \cite{nielson+1994}.  Because the multiplicative correction does not alter the sign of the regression function, this type of correction is particularly well suited for adjusting non-negative regression functions.   \cite{jones++1995} showed that if the true regression function has four continuous derivatives, then the multiplicative bias reduction is operationally equivalent to using an order four kernel.  And while this does remove the bias, it also increases the variance.
\medskip

Many authors have extended the work of \cite{jones++1995}. \cite{glad1998a,glad1998b} proposes to use a parametrically guided local linear smoother and Nadaraya-Watson smoother by starting with a parametric pilot. This approach is extended to a more general framework which includes both multiplicative and additive bias correction by \cite{mamiul2008} (see also \cite{misuul2010} for an extension to time series conditional variance estimation). For multiplicative bias correction in density estimation, we refer the reader to the recent works of \cite{guhanisc2009}, \cite{hasc2007} and \cite{hiru2009}.
\medskip

Although the bias-variance tradeoff  for nonparametric smoothers is always present in finite samples, it is possible to construct smoothers whose \textit{asymptotic bias} converges to zero while keeping the same asymptotic variance. \cite{hengartner+2008}  has exhibited a nonparametric density estimator based on multiplicative bias correction with that property, and have shown in simulations that their estimator also enjoyed good finite sample properties. \cite{burr2010} present such an estimator for nonparametric regression to estimate energy spectra. They illustrate the benefits of this approach on real and simulated spectra. The goal of this paper is to study the asymptotic properties of this estimator. It is worth pointing out that these properties have already been studied by \cite{nielson+1994} for fixed design and further by \cite{jones++1995}. We emphasize that there are two major differences between our work and that of \cite{jones++1995}.
\begin{itemize}
\item First, we do not add regularity assumptions on the target regression function. In particular, we do not assume that the regression function has four continuous derivatives as in \cite{jones++1995}.
\item Second, we show that the multiplicative bias reduction procedure performs a bias reduction with no cost to the asymptotic variance. It is exactly the same as the asymptotic variance of the local linear estimate.
\end{itemize}
We provide another asymptotic behavior under less restrictive assumptions than in \cite{jones++1995}. Moreover our results and proofs are completely different from the above referenced works.
\medskip

This paper is organized as follows. Section \ref{section:preliminaries} introduces the notation and defines the estimator. Section \ref{section:main} gives the asymptotic behavior of the proposed estimator. A brief simulation study on finite sample comparison is presented in Section \ref{section:simulation}. The interested reader is referred to the Appendix where we have gathered the technical proofs.

\section{Preliminaries\label{section:preliminaries}} 

\subsection{Notation.}
Let $(X_1,Y_1),\ldots,(X_n,Y_n)$ be $n$ independent copies of the pair of random variables $(X,Y)$ with values in ${\mathbb R}\times{\mathbb R}$. We suppose that the explanatory variable $X$ has probability density $f$ and model the dependence of the univariate response variable $Y$ to the explanatory variable $X$ through the nonparametric regression model
\begin{equation} \label{eq:model}
Y = m(X) + \varepsilon.
\end{equation}
We assume that the regression function $m(\cdot)$ is  smooth and that the disturbance  $\varepsilon$ is a mean zero random variable with finite variance $\sigma^2$ that is independent of the covariate $X$. Consider the linear smoothers for the regression function $m(x)$ which we can write as
$$\hat m(x) = \sum_{j=1}^n \omega_j(x;h) Y_j,$$
where the weight function $\omega_j(x;h)$ depends on a tuning parameter $h$,  that we think of as the bandwidth.
\medskip

If the weight functions are such that $\sum_{j=1}^n \omega_j(x;h) = 1 \mbox{ and } \sum_{j=1}^n \omega_j(x;h)^2 = (nh)^{-1} \tau^2 ,$ and if the disturbances satisfy the Lindberg-Feller condition, then the linear smoother obeys the central limit theorem
\begin{equation} \label{eq:clt}
\sqrt{nh} \left ( \hat m(x) - \sum_{j=1}^n w_j(x;h)m(X_j) \right ) \longrightarrow
{\mathcal N}(0,\tau^2).
\end{equation}
We can use (\ref{eq:clt}) to construct asymptotic pointwise confidence intervals for the unknown regression function $m(x)$.  But unless the limit of the scaled bias
$$ b(x) = \lim_{n \longrightarrow \infty}\sqrt{nh} \left ( \sum_{j=1}^n w_j(x;h)m(X_j) - m(x) \right ),$$
which we call the asymptotic bias, is  zero,  the confidence interval 
$$\left[\hat m(x) - Z_{1-\alpha/2} \sqrt{nh} \tau, \hat m(x) + Z_{1-\alpha/2}\sqrt{nh} \tau\right]$$
will not cover asymptotically the true regression function $m(x)$ at the nominal $1-\alpha$ level. The construction of valid pointwise $1-\alpha$ confidence intervals for regression smoothers is another motivation for developing estimators with zero asymptotic bias.

\subsection{Multiplicative bias reduction}
Here we present a framework for multiplicative bias reduction. Given a pilot smoother
$$\tilde m_n(x) =\sum_{j=1}^n \omega_{j}(x;h_0) Y_j,$$
the ratio 
$$
V_j = \frac{Y_j}{\tilde m_n(X_j)}
$$
is a noisy estimate of $m(X_j)/\tilde{m}_n(X_j)$, the inverse relative estimation error of the smoother $\tilde m_n$ at each of the observations. Smoothing $V_j$ by
%The inverse relative estimation error can be estimated by smoothing $V_j$, that is,
$$
\widehat \alpha_n(x) = \sum_{j=1}^n \omega_{j}(x;h_1) V_j
$$
yields an estimate for the inverse of the relative estimation error which can be used as a multiplicative correction of the pilot smoother. This leads to the (nonlinear) smoother 
\begin{equation} 
\label{eq:corrected}
\widehat m_n(x) = \widehat \alpha_n(x) \tilde m_n(x).
\end{equation}
%which has smaller bias than $\tilde m(x)$ while maintaining essentially the same variance. 
The estimator (\ref{eq:corrected}) was studied for fixed design by \cite{nielson+1994} and further studied by \cite{jones++1995}.  In both cases, they assumed that the regression function had four continuous derivatives, and show an improvement in the convergence rate of the corrected estimator. \cite{glad1998a,glad1998b} proposed to use a parametrically guided local linear smoother and Nadaraya-Watson smoother by starting with a parametric pilot. She shows that the resulting estimates improve on the local polynomial estimate as soon as the pilot captures some of the features of the regression function.

\section{Theoretical Analysis of Multiplicative 
Bias Reduction\label{section:main}}

In this section, we show that the multiplicative smoother has smaller bias with essentially no cost to the variance, assuming only two derivatives of the regression function. While the derivation of our results are for local linear smoothers, the technique used in the proofs can be easily adapted for other linear smoothers, and the conclusions remain essentially unchanged.

\subsection{Assumptions}
\label{sec:assupmtions}
We make the following assumptions:
\begin{enumerate}
\item The regression function is bounded and strictly positive, that is,
$b\geq m(x) \geq a > 0$ for all $x$.
\item The regression function is twice continuously differentiable everywhere.
\item The density of the covariate is strictly positive on the interior of its support in the sense that $f(x) \geq b({\mathcal K}) > 0$ over every compact ${\mathcal K}$ contained in the support of $f$.
\item $\varepsilon$ has finite fourth moments and has a symmetric distribution around zero.
\item Given a symmetric probability density $K(\cdot)$, consider the weights $\omega_j(x;h)$ associated to the local linear smoother. That is, denote by $K_h(\cdot) = K(\cdot/h)/h$ the scaled kernel by the bandwidth $h$ and define for $k=0,1,2,3$ the sums
$$
S_k(x) \equiv S_k(x;h) = \sum_{j=1}^n (X_j-x)^kK_h(X_j-x).
$$
Then
$$
\omega_{j}(x;h) = \frac{S_2(x;h) - (X_j-x)S_1(x;h)}{S_2(x;h)S_0(x;h)-S_1^2(x;h)} \, K_h(X_j-x).
$$
We set
$$
\omega_{0j}(x) = \omega_j(x;h_0) \quad \mbox{and} \quad
\omega_{1j}(x) = \omega_j(x;h_1).
$$
\item The bandwidths $h_0$ and $h_1$ are such that
$$
h_0 \to 0,  \quad h_1 \to 0, \quad
nh_0 \to \infty, \quad 
nh_1^3 \to \infty,\quad \frac{h_1}{h_0} \to0\quad\textrm{as}\quad n\to \infty.
$$
\end{enumerate}

The positivity assumption (assumption 1) on $m(x)$ is classical when we perform a multiplicative bias correction. It allows to avoid that the terms $Y_j/\tilde{m}_n(X_j)$ blows up. Of course, the regression function might cross the $x$-axis. For such a situation, \cite{glad1998a} proposes to shift all response data $Y_i$ a distance $a$, so that the new regression function $m(x)+a$ does not any more intersect with the $x$-axis. Such a method can also be performed here. Assumptions 2--4 are standard to obtain rate of convergence for nonparametric estimators. Assumption 5 means that we conduct the theory for the local linear estimate. The results can be generalized to other linear smoothers. Assumption 6 is not restrictive since it is satisfied for a wide range of values of $h_0$ and $h_1$.

\subsection{A technical aside}

The proof of the main results rests on establishing a stochastic approximation of estimator (\ref{eq:corrected}) in which each term can be directly analyzed.
\begin{prop}
\label{stochastic}
We have
$$\widehat m_n(x)=\mu_n(x) + \sum_{j=1}^n \omega_{1j}(x)A_j(x)+\sum_{j=1}^ n \omega_{1j}(x)B_j(x) + \sum_{j=1}^n \omega_{1j}(x)\xi_j,$$
where $\mu_n(x)$, conditionally on $X_1,\hdots,X_n$ is a deterministic function,
$A_j$, $B_j$ and $\xi_j$ are random variables. Under condition $nh_0 \rightarrow \infty$, the remainder $\xi_j$ converges to 0 in probability and we have
$$\widehat m_n(x)=\mu_n(x) + \sum_{j=1}^n \omega_{1j}(x)A_j(x) + \sum_{j=1}^n \omega_{1j}(x)B_j(x) + \mathrm{O}_P\left(\frac{1}{nh_0}\right).$$
\end{prop}
\begin{remark}
A technical difficulty arises because even though $\xi_j$ may be small in probability, its expectation may not be small.  We resolve this problem by showing that we only needs to modify $\xi_j$ on a set of vanishingly small probability to guarantee that its expectation is also small.
\end{remark}
%{\bf Remark:} A technical difficulty arises because even though $\xi_j$ may be small in probability, its expectation may not be small.  We resolve this problem by showing that we only needs to modify $\xi_j$ on a set of vanishingly small probability to guarantee that its expectation is also small.

\medskip
\noindent
\begin{defi}
Given a sequence of real numbers $a_n$, say that a sequence of random variables $\xi_n = \mathrm{o}_p(a_n)$ if for all fixed $t > 0$, 
$$
\limsup_{n \longrightarrow \infty} \PP[ |\xi_n| > t a_n] = 0.
$$
\end{defi}
\bigskip

We will need the following Lemma.

\begin{lemma}\label{proba}
If $\xi_n = \mathrm{o}_p(a_n)$, then there exists a sequence of
random variables $\xi^\star_n$ such that
\[
\limsup_{n \longrightarrow \infty} \PP[ \xi_n^\star \not = \xi_n ] = 0 \quad
\mbox{and} \quad \E[\xi_n^\star] = \mathrm{o}(a_n).
\]
We shall use the following notation 
$$
\Es [ \xi_n ] = \E[ \xi^\star_n ].
$$ 
\end{lemma}

\subsection{Main results}
We deduce from Proposition \ref{stochastic} and Lemma \ref{proba} the following Theorem.
\begin{theo}
\label{main}
Under the assumptions (1)-(6), the estimator $\widehat m_n$ satisfies:
$$
%\label{bias}
\Es\left(\widehat m_n(x)|X_1,\hdots,X_n\right)=\mu_n(x)+\mathrm{O}_p\left(\frac{1}{n\sqrt{h_0h_1}}\right) + \mathrm{O}_p\left(\frac{1}{nh_0}\right)
$$
and
$$
\Vs(\widehat m_n(x)|X_1,\hdots,X_n)=\sigma^2\sum_{j=1}^nw_{1j}^2(x) + \mathrm{O}_p\left(\frac{1}{nh_0}\right)+\mathrm{o}_p\left(\frac{1}{nh_1}\right).
$$
\end{theo}
We deduce from Theorem \ref{main} that if the bandwidth $h_0$ of the pilot estimator converges to zero much slower than $h_1$, then $\widehat m_n$ has exactly the same asymptotic variance as the local linear smoother of the original data with bandwidth $h_1$. However, for finite samples, the two step local linear smoother can have a slightly larger variance depending on the choice of $h_0$. For the bias term, a limited Taylor expansion of $\mu_n(x)$ leads to the following result.
\begin{theo}
\label{main2}
Under the assumptions (1)-(6), the estimator $\widehat m_n$ satisfies:
$$
\Es\left(\widehat m_n(x)|X_1,\hdots,X_n\right)=m(x) +\mathrm{o}_p(h_1^2). 
$$
\end{theo}
It is worth pointing out that we only suppose that the regression function is twice continuously differentiable. We do not add smoothness assumptions. For a study of the local linear estimate in the presence of jumps in the derivative, we refer the reader to \cite{desgij2009}. Combining Theorem \ref{main} and Theorem \ref{main2}, we conclude that the multiplicative adjustment performs a bias reduction on the pilot estimator without increasing the asymptotic variance. The asymptotic behavior of the bandwidths $h_0$ and $h_1$ is constrained by assumption 6. However, it is easily seen that this assumption is satisfied for a large set of values of $h_0$ and $h_1$. For example, the choice $h_1=c_1n^{-1/5}$ and $h_0=c_0n^{-\alpha}$ for $0<\alpha<1/5$ leads to
$$\Es\left(\widehat m_n(x)|X_1,\hdots,X_n\right)-m(x)=\mathrm{o}_p(n^{-2/5})$$
and
$$\Vs(\widehat m_n(x)|X_1,\hdots,X_n)=\mathrm{O}_p\left(n^{-4/5}\right).$$
\begin{remark}
Estimators with bandwidths of order $\mathrm{O}(n^{-\alpha})$ for $0<\alpha<1/5$ are oversmoothing the true regression function, and as a result, they have biases that are of larger order of magnitude than their standard deviations. We conclude that the multiplicative adjustment performs a bias reduction on the pilot estimator.
\end{remark}

\section{Numerical examples\label{section:simulation}}
While the amount of the bias reduction depends on the curvature of the regression function, a decrease is expected (asymptotically) everywhere, and this, at no cost to the variance. The simulation study in this section shows that this asymptotic behavior emerges already at modest sample sizes.

\subsection{Local study}
\begin{figure}[h]
  \centering
  \includegraphics[width=8cm,height=10cm,angle=-90]{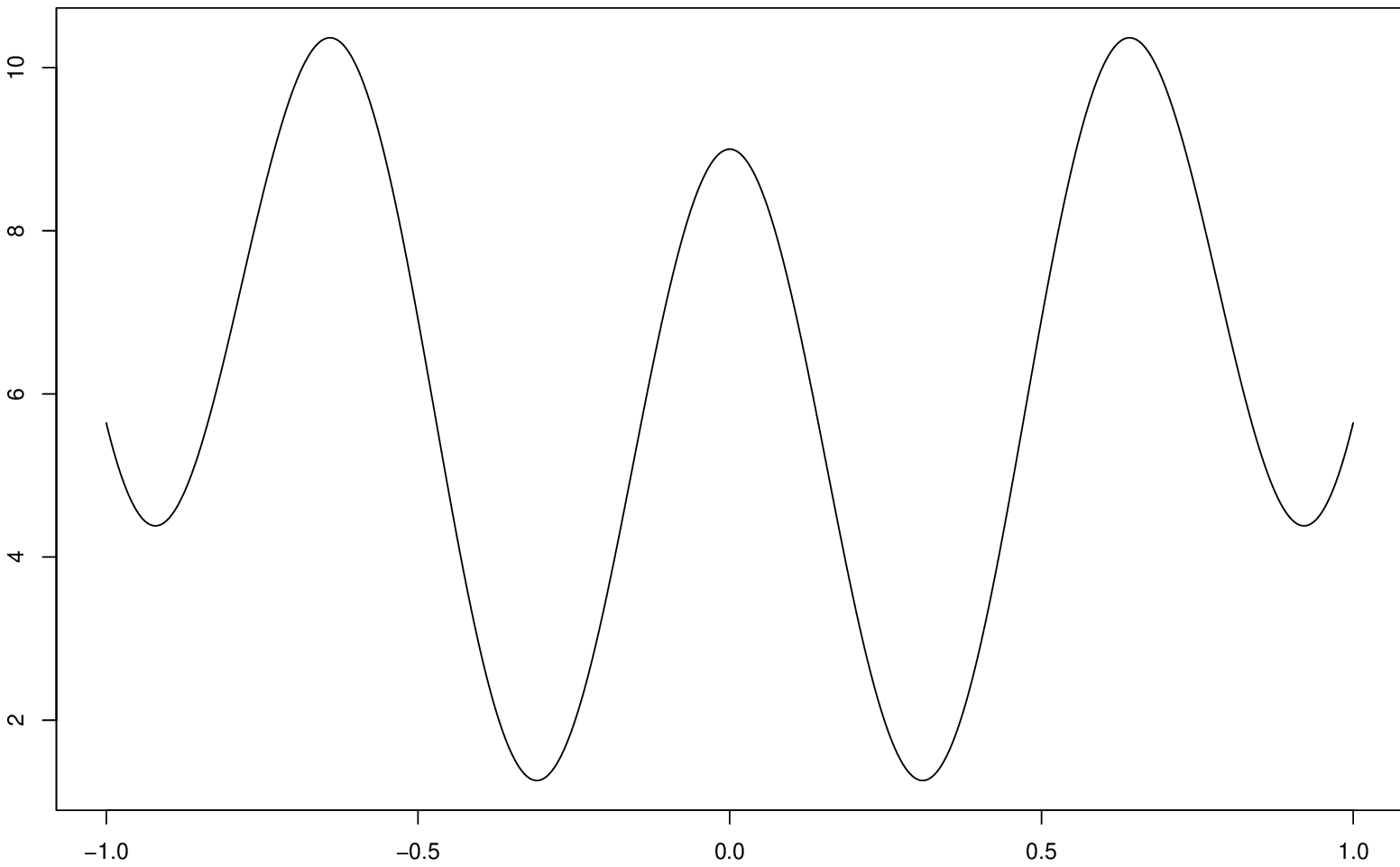}
  \caption{\textsf{The regression function to be estimated.}}
  \label{fig:num_ex}
\end{figure}

To illustrate numerically the possible reduction in the bias and associate increase of the variance achieved by the multiplicative bias correction, consider estimating the regression function
$$m(x)=5+3|x|^{5/2}+x^2+4\cos(10x)$$
at $x=0$ (see Figure \ref{fig:num_ex}).
The local linear smoother tends to under-estimate the regression function at their maximum, and hence, this example will provide a good example. Furthermore, because the second derivative of this regression function is continuous but not differentiable at the origin, the results previously obtained by \cite{nielson+1994} do not apply.
\medskip

The data are simulated according to the model
$$Y_i=m(X_i)+\varepsilon_i,\quad i=1,\hdots,100,$$
where $\varepsilon_i$ are independent $\mathcal N(0,1)$ variables. We first consider the local linear estimate with a Gaussian kernel function and we study its performances over a grid of bandwidths $\mathcal H=[0.005,0.1]$. For the new estimate, the theory recommends to start with an oversmooth pilot estimate. In this regard, we take $h_0=0.1$ and study the performance of the multiplicative bias corrected estimate for $h_1\in\mathcal H_1=[0.005,0.12]$. To explore the stability of our two stages estimator with respect to $h_0$, we also consider the choice $h_0=0.02$. For such a choice, the pilot estimate clearly undersmoothes the regression function.
\medskip

\begin{figure}[H]
  \centering
  \includegraphics[width=13.5cm,height=8cm]{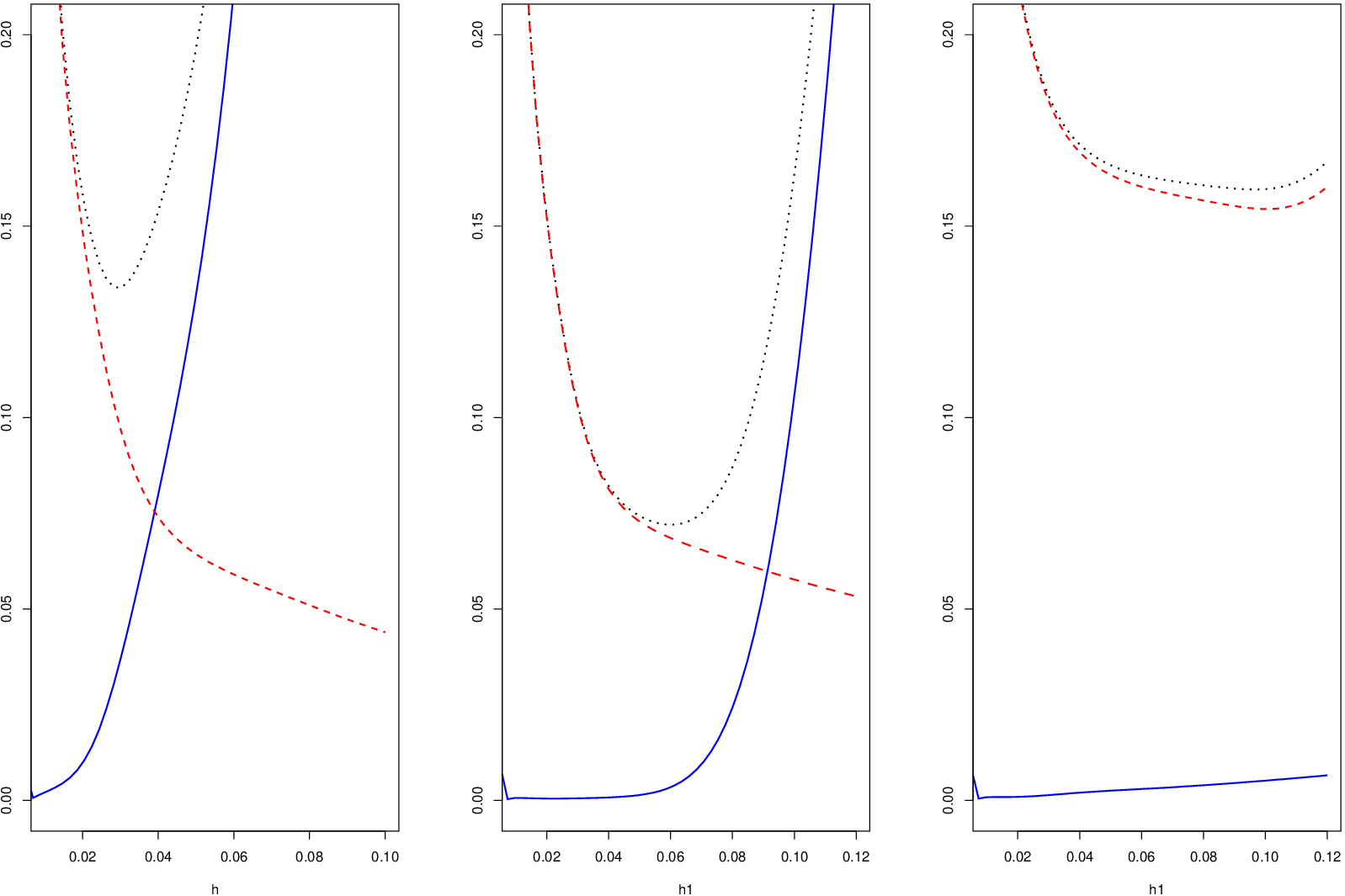}
  \caption{\textsf{Mean square error (dotted line), squared bias (solid line) and variance (dashed line) of the local linear estimate (left) and multiplicative bias corrected estimate with $h_0=0.1$ (center) and $h_0=0.02$ (right) at point $x=0$.}}
  \label{fig:bias_variance}
\end{figure}

Bias and variance of each estimate are calculated at $x=0$. To do this, we compute the value of each estimate at $x=0$ for 200 samples $(X_i,Y_i),i=1,\hdots,100$. The same design $X_i,i=1,\hdots,100$ is used for each sample. It is generated according to a uniform distribution over $[-1,1]$. The bias at point $x=0$ is estimated by subtracting $m(0)$ at the mean value of the estimate at $x=0$ (the mean value is computed over the 200 replications). Similarly we estimate the variance at $x=0$ by the variance of the values of the estimate at this point. Figure \ref{fig:bias_variance} presents squared bias, variance and mean square error of each estimate for different values of bandwidths $h$ for the local linear smoother and $h_1$ for our estimate.

The first conclusion is that the corrected estimate has smaller bias than the local linear estimate provided the pilot estimate oversmoothes the regression function. Small values of $h_0$ clearly undersmooth the regression function, whatever the choice of $h_1$. Moreover, it is worth pointing out that our procedure does not significantly increase the variance. Even if Theorem \ref{main} and Theorem \ref{main2} provide asymptotic results, our simulations show that the asymptotic behavior of our estimate emerges already at modest sample size. Finally, due to the bias reduction, we note that our procedure also reduces the optimal mean square error (see Table \ref{tab:comp_opt_mse}).

\begin{table}[H]
  \centering
  \begin{tabular}{|c||c|c|c|}
\hline
 & MSE & Bias$^2$ & Variance \\
\hline\hline
LLE & 0.134 & 0.038 & 0.096  \\
\hline
MBCE & 0.072 & 0.003 & 0.069 \\
\hline
  \end{tabular}
  \caption{\textsf{Optimal mean square error (MSE) for the local linear estimate (LLE) and the multiplicative bias corrected estimate (MBCE) with $h_0=0.1$ at point $x=0$.}}
  \label{tab:comp_opt_mse}
\end{table}

We conclude our local study with a comparison between our estimate and the estimate proposed by \cite{glad1998a}. To do this, we compute the multiplicative bias corrected estimate using three parametric starts:
\begin{itemize}
\item first the guide is chosen correctly and belong to the true parametric family:
$$\tilde m_n^1(x)=\hat\beta_0+\hat\beta_1|x|^{5/2}+\hat\beta_2 x^2+\hat\beta_3\cos(10x);$$
\item second, we consider a linear parametric guide (which is obviously wrong):
$$\tilde m_n^2(x)=\hat\beta_0+\hat\beta_1 x;$$
\item finally, we use a more reasonable guide, not correct, but that can reflect some a priori idea on the regression function
$$\tilde m_n^3(x)=\hat\beta_0+\hat\beta_1x+\hat\beta_2 x^2+\hdots+\hat\beta_8x^8.$$
\end{itemize}
All the estimates $\hat\beta_j$ stands for the classical least square estimates.
\medskip

The multiplicative bias correction is performed on these parametric starts using the local linear estimate. The performance of the resulting estimates is measured over a grid of bandwidths $\mathcal H_2=[0.005;0.4]$. Bias and variance of each estimate are still estimated at $x=0$. We keep the same setting as above and all the results are averaged over the same 200 replications. We display in Table \ref{tab:comp_opt_mse_param} the optimal MSE calculated over the grid $\mathcal H_2$.

\begin{table}[H]
  \centering
  \begin{tabular}{|c||c|c|c|}
\hline
 & MSE & Bias$^2$ & Variance \\
\hline\hline
start $\tilde m_n^1$ & 0.060 & 0.000 & 0.060 \\
\hline
start $\tilde m_n^2$ & 0.134 & 0.038 & 0.096 \\
\hline
start $\tilde m_n^3$ & 0.095 & 0.021 & 0.074 \\
\hline
  \end{tabular}
  \caption{\textsf{Optimal mean square error for the multiplicative bias corrected estimates with parametric starts $\tilde m_n^j$, $j=1,2,3$.}}
  \label{tab:comp_opt_mse_param}
\end{table}
As expected, we first observe that the performance clearly depends on the choice of the parametric start. Table \ref{tab:comp_opt_mse} and table \ref{tab:comp_opt_mse_param} show that (in term of MSE) the estimate studied in this paper is better than the corrected estimated with parametric start $\tilde m_n^2$ and $\tilde m_n^3$. Unsurprisingly, the best performance are obtained with the parametric guide $\tilde m_n^1$ (which belongs to the true model). In practice, when one has no or few a priori information on the target regression function, the method proposed in the present paper is preferable.

\subsection{Global study}

This paper does not conduct any theory to select the two bandwidths $h_0$ and $h_1$ in an optimal way. If automatic procedures are needed, they can be obtained by adjusting traditional automatic selection procedures for the classical nonparametric estimators (see \cite{burr2010}). In this part, we propose to use leave-one-out cross validation to choose both $h_0$ and $h_1$. We then compare the performance of the selected estimate with the local polynomial estimate in term of integrated square error.
\medskip

\begin{figure}[h]
  \centering
  \includegraphics[width=13.5cm,height=8cm]{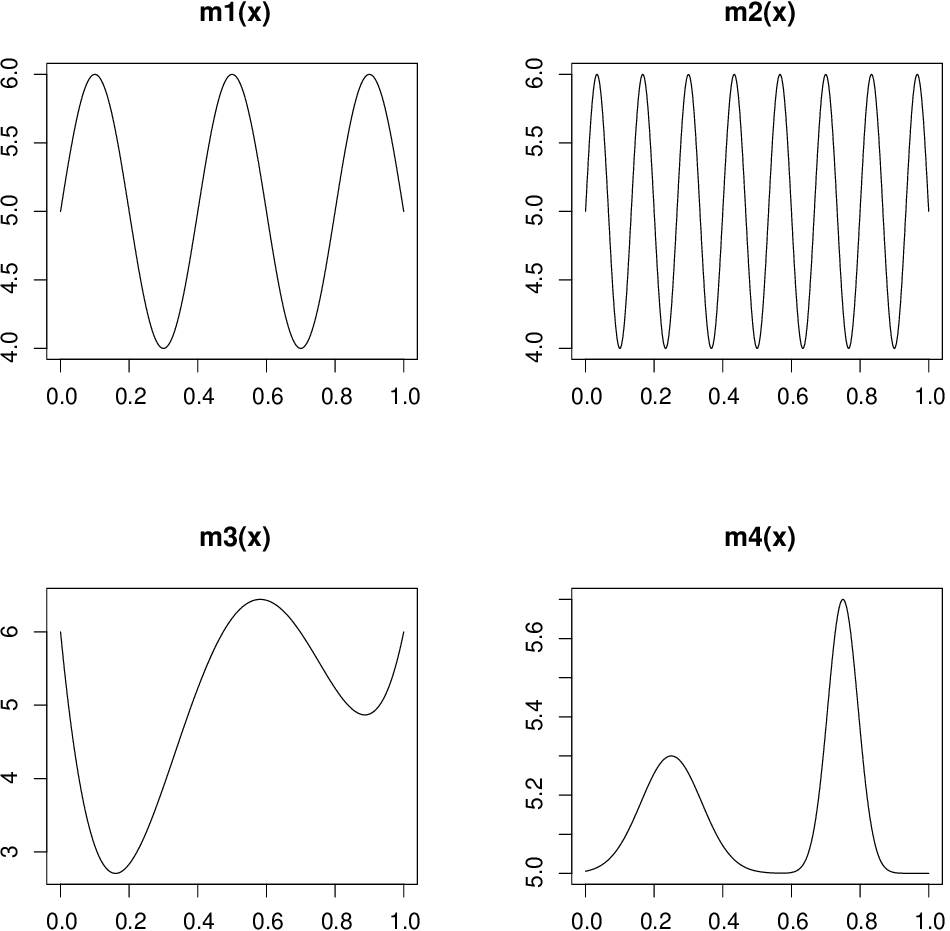}
  \caption{\textsf{Regression functions to be estimated.}}
  \label{fig:4fonc_tests}
\end{figure}

\cite{hurvich++1998} report a comprehensive numerical study that compares standard smoothing methods on various test functions. Here, we take the same setting to compare the local linear estimate with its multiplicative bias corrected smoother. In each of the examples, we take the Gaussian kernel $K(x)=\exp(-x^2/2)/\sqrt{2\pi}$. We use the following regression functions (see Figure \ref{fig:4fonc_tests}):
\medskip

\begin{tabular}{ccl}
\hspace*{5mm}&(1)& $m_1(x)=\sin(5\pi x)$\\
&(2)& $m_2(x)=\sin(15\pi x)$\\
&(3)& $m_3(x)=1-48x+218x^2-315x^3+145x^4$\\
&(4)& $m_4(x)=0.3\exp{[-64(x-.25)^2]}+0.7\exp{[-256(x-.75)^2]}$\\
%&(5)& $m_5(x)=10\exp{(-10x)}$\\
%&(6)& $m_6(x)=\exp{(x-1/3)}\{x<1/3\}+\exp{-2(x-1/3)}\{x\geq1/3\}$
\end{tabular}
\smallskip

and we take a Gaussian error distribution with standard deviation $\sigma_j=0.25\|m_j\|_2$ for $j=1,\hdots,4$.
\medskip

We use a cross validation device to select both $h_0$ and $h_1$. This selection procedure involves solving minimization problem that necessitate a search over a finite grid $\H$ of bandwidths $h_0$ and $h_1$. Formally, given $\H$, we choose $\hat h_0$ and $\hat h_1$ such as
$$
(\hat h_0,\hat h_1)= \argmin_{(h_0,h_1) \in \mathcal \H \times \H}  
\frac{1}{n}\sum_{i=1}^n(Y_i-\widehat m_n^i(X_i))^2.
$$
Here $\widehat m_n^i$ stands for the corrected local polynomial estimate after deleted the $i$th observation. To assess the quality of the selected estimate, we compare its performances with the local polynomial estimate for which the bandwidth is again selected by leave-one-out cross validation. The performance of an estimator $\widehat m$ is measured by the integrated square error 
$$ISE(\widehat m)=\int_0^1\left(m(x)-\widehat m(x)\right)^2\,\textrm{d}x,$$
and to avoid the boundary effects, the design $X_1,\hdots,X_n$ is generated according to a uniform distribution over $[-0.2,1.2]$. 
\bigskip

\begin{figure}[H]
  \centering
  \includegraphics[width=10cm,height=12cm,angle=-90]{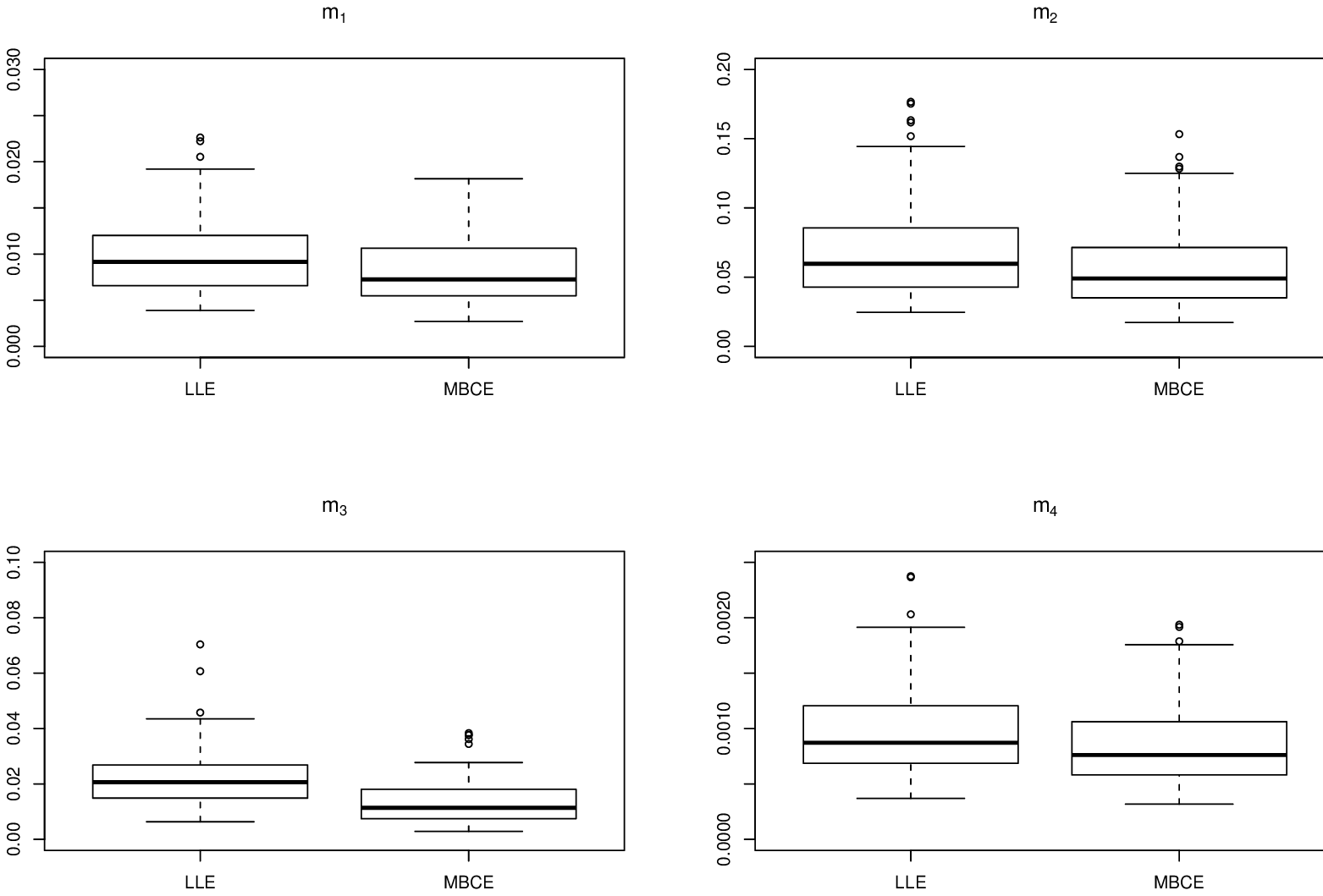}
  \caption{\textsf{Boxplot of the integrated square error over the 100 replications.}}
  \label{fig:boxplot_ISE_VC}
\end{figure}

Table \ref{tab:res_global} presents the median over 100 replications of
\begin{itemize}
\item the selected bandwidths;
\item the integrated square error;
\item the integrated square error of the local linear estimate divided by  the integrated square error of the corrected estimate ($R_{ISE}$).
\end{itemize}
Figure \ref{fig:boxplot_ISE_VC} displays the boxplots of the integrated square error for each estimate.

\begin{table}[H]
  \centering
  \begin{tabular}{|c||c|c||c|c|c||c|}
\hline
 & \multicolumn{2}{|c||}{LLE} & \multicolumn{3}{|c||}{MBCE} & \\
\hline\hline
 & $h$ & ISE {\footnotesize ($\times 100$)} & $h_0$ & $h_1$ & ISE {\footnotesize ($\times 100$)} & $R_{ISE}$ \\
\hline\hline
$m_1$ & 0.023 & 0.920 & 0.041 & 0.032 & 0.727 & 1.191 \\
$m_2$ & 0.011 & 5.967 & 0.027 & 0.012 & 4.968 & 1.205 \\
$m_3$ & 0.029 & 2.063 & 0.071 & 0.054 & 1.139 & 1.648 \\
$m_4$ & 0.018 & 0.087 & 0.033 & 0.023 & 0.076 & 1.147 \\
\hline
  \end{tabular}
  \caption{\textsf{Median over 100 replications of the selected bandwidths and of the integrated square error of the selected estimates. LLE and MBCE stands for local linear estimate and multiplicative bias corrected estimate.}}
  \label{tab:res_global}
\end{table}

We obtain significant ISE reduction for the four models. As predicted by Theorem \ref{main}, the data-driven procedure selects $h_0$ bigger than $h$: the pilot estimate is oversmoothing the true regression function. Of course, selecting both $h_0$ and $h_1$ is time consuming and can appear as the price to be paid to improve the local linear smoother.
\medskip

\begin{figure}[H]
  \centering
  \includegraphics[height=12cm,width=8cm,angle=270]{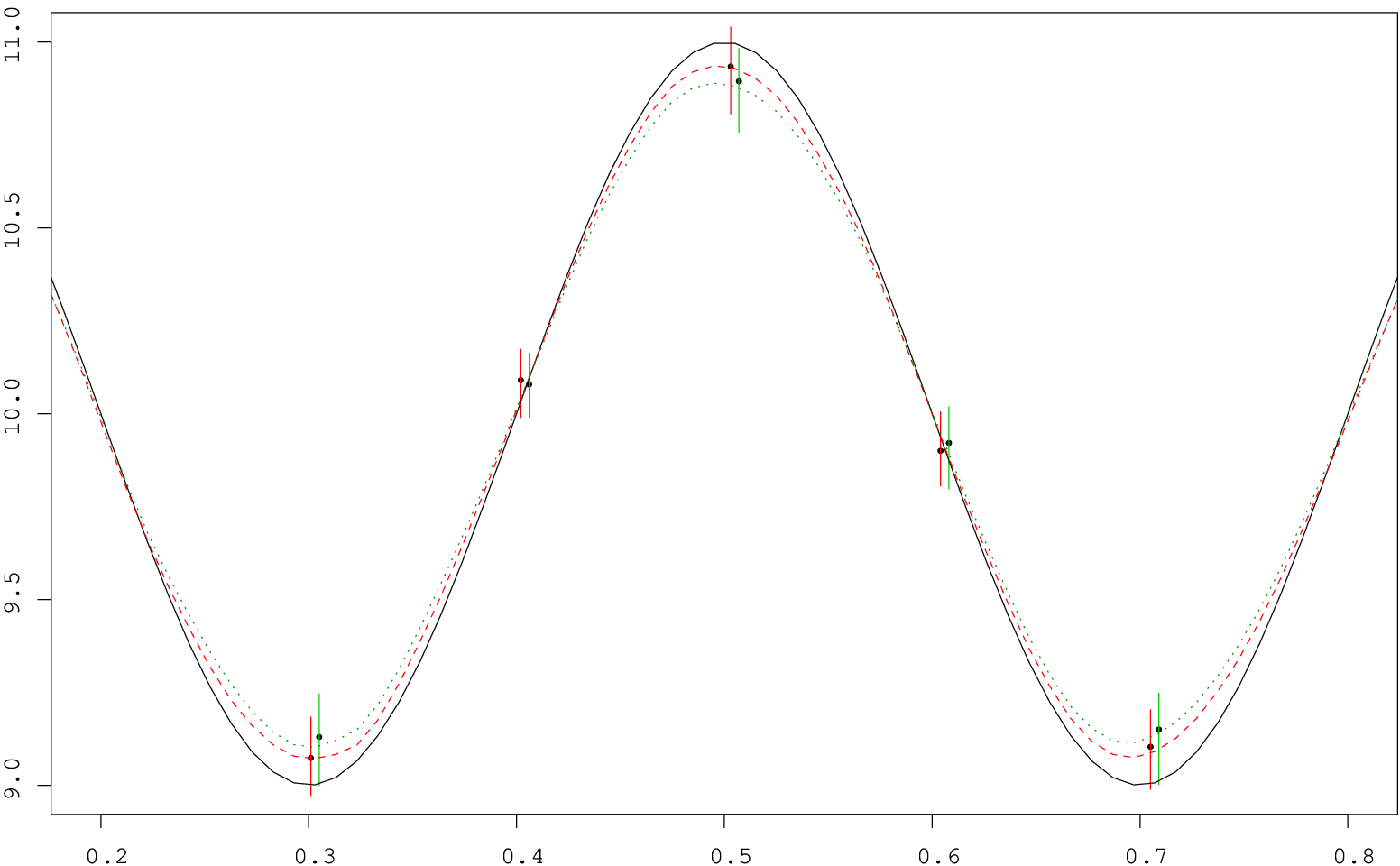}
  \caption{\textsf{The solid curve represents the true regression function, our estimator is in dashed line and local linear smoother is dotted.}}
  \label{fig:biais_var_m2}
\end{figure}

Figure \ref{fig:biais_var_m2} presents, for the regression function $m_1$ with $n=100$ and 100 iterations, different estimators on a grid of points. In lines is the true regression function which is unknown. For every point on a fixed grid, we plot, side by side, the mean over 100 replications of our estimator at that point (left side) and on the right side of that point the mean over 100 replications of the local polynomial estimator. Leave-one-out cross validation is applied to select the bandwidths $h_0$ and $h_1$ for our estimator and the bandwidth $h$ for the local polynomial estimator. We add also the interquartile interval in order to see the fluctuations of the different estimators.
In this example, our estimator reduces the bias by increasing the peak and decreasing the valleys. Moreover, the interquartile intervals look similar for both estimator, as predicted by the theory.

\section{Proofs}
This section is devoted the technical proofs. 
\subsection{Proof of Proposition \ref{stochastic}}
Write the bias corrected estimator 
$$\widehat m_n(x)=\sum_{j=1}^n \omega_{1j}(x)\frac{\tilde m_n(x)}{\tilde m_n(X_j)}Y_j=\sum_{j=1}^n \omega_{1j}(x) R_j(x) Y_j,
$$
and let us approximate the quantity $R_j(x)$. Define
$$
\bar m_n(x)=\sum_{j=1}^n \omega_{0j}(x) m(X_j)= \E\left(\tilde m_n(x)|X_1,\hdots,X_n\right),
$$
and observe that
\begin{align*}
R_j(x)=& \frac{\tilde m_n(x)}{\tilde m_n(X_j)}\\
=& \frac{\bar m_n(x)}{\bar m_n(X_j)} \times 
\left ( 1 + \frac{\tilde m_n(x)-\bar m_n(x)}{\bar m_n(x)}\right ) \times
\left (1 +  \frac{\tilde m_n(X_j)-\bar m_n(X_j)}{\bar m_n(X_j)}\right )^{-1}\\
=& \frac{\bar m_n(x)}{\bar m_n(X_j)} \times 
[1+\Delta_n(x)] \times \frac{1}{1+\Delta_n(X_j)},
\end{align*}
where
$$
\Delta_n(x)= \frac{\tilde m_n(x)-\bar m_n(x)}{\bar m_n(x)} =
\frac{\sum_{l \leq n} \omega_{0l}(x) \varepsilon_l}{\sum_{l \leq n} \omega_{0l}(x) m(X_l)}.
$$
Write now $R_j(x)$ as 
$$
R_j(x)=\frac{\bar m_n(x)}{\bar m_n(X_j)}
\left[1+\Delta_n(x)-\Delta_n(X_j)+r_j(x,X_j)\right]
$$
where $r_j(x,X_j)$ is a random variable converging to 0 to be define latter on. Given the last expression and model \eqref{eq:model},  estimator (\ref{eq:corrected}) could be written as 
\begin{align*}
\widehat m_n(x)=& \sum_{j=1}^n \omega_{1j}(x)R_j(x)Y_j \\
=&  \sum_{j=1}^n \omega_{1j}(x)\frac{\bar m_n(x)}{\bar m_n(X_j)}m(X_j)\\
\quad& +  \sum_{j=1}^n \omega_{1j}(x)\frac{\bar m_n(x)}{\bar m_n(X_j)}\left[
\varepsilon_j + m(X_j)\left(\Delta_n(x)-\Delta_n(X_j)\right)\right]\\
\quad& +  \sum_{j=1}^n \omega_{1j}(x)\frac{\bar m_n(x)}{\bar m_n(X_j)}\left(
\Delta_n(x)-\Delta_n(X_j)\right) \varepsilon_j\\
\quad& +  \sum_{j=1}^n \omega_{1j}(x)\frac{\bar m_n(x)}{\bar m_n(X_j)}r_j(x,X_j)Y_j\\
=& \mu_n(x) + \sum_{j=1}^n \omega_{1j}(x)A_j(x)+\sum_{j=1}^n\omega_{1j}(x)B_j(x)+\sum_{j=1}^n\omega_{1j}(x)\xi_j.
\end{align*}
which is the first part of the proposition. Under assumption set forth in Section \ref{sec:assupmtions}, the pilot smoother $\tilde m_n$ converges to the true regression function $m(x)$. \cite{bickel1973} show that this convergence is uniform over compact sets ${\mathcal K}$ contained in the support of the density of the covariate $X$. As a result 
$$
\sup_{x \in {\mathcal K}} | \tilde m_n(x) - \bar m_n(x) | \leq \frac{1}{2}.
$$
So a limited expansion of $(1+u)^{-1}$ yields for $x \in \mathcal{K}$
$$R_j(x)=\frac{\bar m_n(x)}{\bar m_n(X_j)} \left[ 1 + \Delta_n(x) - \Delta_n(X_j)+ \mathrm{O}_p\left(|\Delta_n(x) \Delta_n(X_j)| + \Delta_n^2(X_j)\right) \right],
$$
thus 
$$\xi_j=\mathrm{O}_p\left(|\Delta_n(x) \Delta_n(X_j)| + \Delta_n^2(X_j)\right).$$
Under the stated regularity assumptions, we deduce that
$$\xi_j=\mathrm{O}_p\left(\frac{1}{nh_0}\right).$$
leading to the announced result. Proposition (\ref{stochastic}) is proved.

\subsection{Proof of lemma (\ref{proba})}
By definition
$$\limsup_{n \longrightarrow \infty} \PP[ |\xi_n| > t a_n] = 0$$ 
for all $t>0$,  so that a triangular array argument shows that there exists
an increasing sequence $m=m(k)$ such that
$${\mathbb P} \left [|\xi_n| > \frac{a_n}{k} \right ] \leq \frac{1}{k} \qquad \mbox{for all } n \geq m(k).$$ 
For $ m(k) \leq n \leq m(k+1)-1$, define
$$
\xi^\star_n = \left \{ \begin{array}{ll} \xi_n & \mbox{ if } |\xi_n| < k^{-1} a_n \\
0 & \mbox{otherwise.}
\end{array} \right .
$$
It follows from the construction of $\xi^\star_n$ that for $n \in (m(k),m(k+1)-1)$,
$${\mathbb P}[\xi_n \not = \xi_n^*] = {\mathbb P}[|\xi_n| > k^{-1} a_n] \leq \frac{1}{k},
$$
which converges to zero as $n$ goes to infinity.  Finally set $k(n)=\sup\{k:m(k)\leq n\}$, we obtain
$$\E[|\xi|_n^\star] \leq \frac{a_n}{k(n)} = \mathrm{o}(a_n).$$

\subsection{Proof of Theorem (\ref{main})}
Recall that 
$$
\widehat m_n(x)=\mu_n(x) + \sum_{j=1}^n \omega_{1j}(x)A_j(x) 
+ \sum_{j=1}^n \omega_{1j}(x)B_j(x) + \mathrm{O}_P\left(\frac{1}{nh_0}\right).
$$
Focus on the conditional bias, we get 
\begin{eqnarray*}
\E(\mu_n(x)| X_1,\ldots,X_n) &=& \mu_n(x)\\
\E(A_j(x) | X_1,\ldots,X_n) &=& 0\\
\E(B_j(x) | X_1,\ldots,X_n) &=&\frac{\bar m_n(x)}{\bar m_n(X_j)}\sigma^2 
\Big(\frac{\omega_{0j}(x)}{\bar m_n(x)} -  \frac{\omega_{0j}(X_j)}{\bar m_n(X_j)}\Big).
\end{eqnarray*}
Since
$$
\sum_{j=1}^n \omega_{1j}(x)\omega_{0j}(x) \leq  \sqrt{\sum_{j=1}^n \omega_{1j}(x)^2}\sqrt{\sum_{j=1}^n \omega_{0j}(x)^2}= \mathrm{O}_p\left(\frac{1}{n\sqrt{h_0h_1}}\right),
$$
we deduce that
$$\E\left(\sum_{j=1}^n \omega_{1j}(x)B_j(x)\Big|X_1,\hdots,X_n\right)=\mathrm{O}_p\left(\frac{1}{n\sqrt{h_0h_1}}\right).$$
This proves the first part of the Theorem.
\bigskip

For the conditional variance, we use the following expansion of the two stages estimator
$$
\widehat m_n(x)=\sum_{j=1}^n \omega_{1j}(x)\frac{\bar m_n(x)}{\bar m_n(X_j)}Y_j\left(1+ \left[\Delta_n(x)-\Delta_n(X_j) \right]\right) + \mathrm{O}_p\left(\frac{1}{nh_0}\right).
$$
Using the fact that the residuals have four finite moments and have a symmetric distribution around 0, a moment's thought shows that
$$\V(Y_j\left[\Delta_n(x)-\Delta_n(X_j) \right]| X_1,\ldots,X_n )=\mathrm{O}_p\left(\frac{1}{nh_0}\right)$$
and
$$\C(Y_j,Y_j\left[\Delta_n(x)-\Delta_n(X_j)\right]| X_1,\ldots,X_n )=\mathrm{O}_p\left( \frac{1}{nh_0} \right).$$
Hence
\begin{align*}
\Vs(\widehat m_n(x)| X_1,\ldots,X_n ) & =\V\left(\sum_{j=1}^n \omega_{1j}(x)\frac{\bar m_n(x)}{\bar m_n(X_j)}Y_j\Big| X_1,\hdots,X_n\right) \\
& \hspace{6cm}+\mathrm{O}_p\left(\frac{1}{nh_0}\right).
\end{align*}
Observe that the first term on the right hand side of this equality can be 
seen as the variance of the two stages estimator with a deterministic pilot 
estimator. It follows from \cite{glad1998a} that
$$
\V\left(\sum_{j=1}^n \omega_{1j}(x)\frac{\bar m_n(x)}{\bar m_n(X_j)}Y_j\Big|X_1,\hdots,X_n\right)=\sigma^2\sum_{j=1}^n\omega_{1j}^2(x)+\mathrm{o}_p\left(\frac{1}{nh_1}\right),
$$
which proves the theorem.

\subsection{Proof of theorem (\ref{main2})}
Recall that
$$\mu_n(x)=\sum_{j \leq n} \omega_{1j}(x)\frac{\bar m_n(x)}{\bar m_n(X_j)} m(X_j).$$

We consider  the limited Taylor expansion of the ratio 
$$\frac{m(X_j)}{\bar m_n(X_j)} = \frac{m(x)}{\bar m_n(x)}+(X_j-x) \left (   \frac{m(x)}{\bar m_n(x)} \right )^\prime +\frac{1}{2} (X_j-x)^2 \left (   \frac{m(x)}{\bar m_n(x)}\right )^{\prime \prime} (1+\mathrm{o}_p(1)),$$
then
\begin{eqnarray*}
\mu_n(x) & = & \bar m_n(x)   \left \{ \frac{m(x)}{\bar m_n(x)}
\sum_{j=1}^n \omega_{1j}(x)  
 + \left (   \frac{m(x)}{\bar m_n(x)} \right )^\prime
 \sum_{j=1}^n (X_j-x) \omega_{1j}(x) \right . \\
&& \left .\qquad + \frac{1}{2}\left (   \frac{m(x)}{\bar m_n(x)} 
\right )^{\prime \prime}  \sum_{j=1}^n (X_j-x)^2 \omega_{1j}(x) (1+\mathrm{o}_p(1)) \right\}.
\end{eqnarray*}
It is easy to verify that
\begin{eqnarray*}
\Sigma_0(x;h_1) &=& \sum_{j=1}^n \omega_{1j}(x) = 1, \\
\Sigma_1(x;h_1) &=&  \sum_{j=1}^n (X_j-x) \omega_{1j}(x) = 0 \\
\Sigma_2(x;h_1) &=& \sum_{j=1}^n (X_j-x)^2 \omega_{1j}(x) = 
\frac{S_2^2(x;h_1)-S_3(x;h_1)S_1(x;h_1)}{S_2(x;h_1)S_0(x;h_1)-S_1^2(x;h_1)}.
\end{eqnarray*}
For random designs, we can further approximate (see, e.g., \cite{wand+1995})
$$S_k(x,h_1)=\left\{
  \begin{array}{ll}
h^k\sigma^k_Kf(x)+\mathrm{o}_p(h^k) & \textrm{for }k\textrm{ even} \\
h^{k+1}\sigma^{k+1}_Kf'(x)+\mathrm{o}_p(h^{k+1}) & \textrm{for }k\textrm{ odd,} \\
  \end{array}\right.
$$
where $\sigma^k_K=\int u^kK(u)\,\mathrm{d}u.$ Therefore
\begin{eqnarray*}
\Sigma_2(x;h_1) &=& h_1^2 \int u^2 K(u)\,\mathrm{d}u + \mathrm{o}_p(h_1^2)\\
&=& \sigma_K^2 h_1^2 + \mathrm{o}_p(h_1^2),
\end{eqnarray*}
so that we can write $\mu_n(x)$ as
\begin{align*}
\mu_n(x) = & \bar m_n(x) \left \{ \frac{m(x)}{\bar m_n(x)} +
\frac{\sigma_K^2 h_1^2}{2}\left (   \frac{m(x)}{\bar m_n(x)} 
\right )^{\prime \prime} + \mathrm{o}_p(h_1^2) \right \} \\
= &  m(x)+\frac{\sigma_K^2 h_1^2}{2} \bar m_n(x)   \left ( \frac{m(x)}{\bar m_n(x)} \right )^{\prime \prime} + \mathrm{o}_p(h_1^2).
\end{align*}
Moreover
\begin{eqnarray*}
\left ( \frac{m(x)}{\bar m_n(x)} \right )^{\prime \prime} &=&
\frac{\bar m_n^2(x) m^{\prime \prime}(x)}{\bar m_n^3(x)} 
- 2\frac{\bar m_n(x) \bar m_n^\prime(x) m^\prime(x)}{\bar m_n^3(x)} \\
&& \qquad
- \frac{m(x) \bar m_n(x) \bar m_n^{\prime \prime}(x)}{\bar m_n^3(x)} 
+ 2 \frac{m(x) (\bar m_n^\prime(x))^2}{\bar m_n^3(x)}
\end{eqnarray*}
and applying the usual approximations, we conclude that
$$
\left ( \frac{m(x)}{\bar m_n(x)} \right )^{\prime \prime} = \mathrm{o}_p(1).
$$
Putting all pieces together, we obtain
$$
\Es(\widehat m_n(x)|X_1,\hdots,X_n) - m(x) = \mathrm{o}_p(h_1^2)  + \mathrm{O}_p\left(\frac{1}{n\sqrt{h_0h_1}}\right)+ \mathrm{O}_p\left(\frac{1}{nh_0}\right).
$$
Since
$$nh_1^3 \longrightarrow \infty \quad\textrm{and}\quad \frac{h_1}{h_0}\longrightarrow 0,$$
we conclude that the bias is of order $\mathrm{o}_p(h_1^2)$.

%\paragraph{Acknowledgments.}  The authors greatly thank the Associate Editor and two referees for valuable comments and insightful suggestions.
%%%%%%%%%%%%%%%%%%%%%%%

\bibliographystyle{abbrvnat} 
\bibliography{./biblio}

\end{document}